\documentclass{amsart}
\usepackage{amsmath,amsthm,amssymb,euscript,multirow,bigstrut}
\usepackage[matrix,arrow]{xy}

\newcommand{\C}{\mathbb{C}}

\newcommand{\No}{\operatorname{N}}
\newcommand{\Z}{\mathbb{Z}}
\newcommand{\Q}{\mathbb{Q}}

\newcommand{\R}{\mathbb{R}}

\renewcommand{\a}{{\mathfrak{a}}}

\newcommand{\p}{{\mathfrak{p}}}
\renewcommand{\P}{{\mathfrak{P}}}
\newcommand{\q}{{\mathfrak{q}}}

\newcommand{\E}{\mathcal{E}}
\newcommand{\A}{\mathcal{A}}

\renewcommand{\Re}{\operatorname{Re}}

\newcommand{\category}[1]{\mbox{\boldmath $\mathsf{{#1}}$}}

\newcommand{\conj}[1]{{\overline{{#1}}}}

\newlength{\algcwidth} \newlength{\algcheight}

\newcommand{\jacobi}[2]{{\genfrac{(}{)}{}{}{#1}{#2}}}
\renewcommand{\H}{\category{H}}
\renewcommand{\O}{\mathcal{O}}

\newcommand{\<}{\langle}
\renewcommand{\>}{\rangle}

\renewcommand{\div}{\mid}
\newcommand{\ndiv}{\nmid}

\newcommand{\Aut}{\operatorname{Aut}}

\newcommand{\End}{\operatorname{End}}

\newcommand{\Cl}{{\operatorname{Cl}}}

\newtheorem{theorem}{Theorem}[section]
\newtheorem{proposition}[theorem]{Proposition}
\newtheorem{lemma}[theorem]{Lemma}

\theoremstyle{definition}
\newtheorem{definition}[theorem]{Definition}
\begin{document}

\title{Supersingular primes for points on $X_0(p)/w_p$}
\author{David Jao}
\date{March 16, 2004}
\address{Microsoft Research \\ 1 Microsoft Way \\ Redmond, WA 98052}
\email{davidjao@microsoft.com}
\begin{abstract}
For small odd primes $p$, we prove that most of the rational points on
the modular curve $X_0(p)/w_p$ parametrize pairs of elliptic curves
having infinitely many supersingular primes.  This result extends the
class of elliptic curves for which the infinitude of supersingular
primes is known.  We give concrete examples illustrating how these
techniques can be explicitly used to construct supersingular primes
for such elliptic curves.  Finally, we discuss generalizations to
points defined over larger number fields and indicate the types
of obstructions that arise for higher level modular curves.
\end{abstract}

\maketitle

\section{Introduction}
Let $E$ be an elliptic curve defined over a number field.  It is
conjectured that $E$ has infinitely many prime ideals of supersingular
reduction. For curves $E$ with complex multiplication, a classical
result of Deuring~\cite{deuring} states that the supersingular primes
have density $1/2$. More recently, Elkies proved that $E$ always has
infinitely many supersingular primes whenever it is defined over a
real number field~\cite{e2}, or when the absolute norm of $j(E) -
1728$ has a prime factor congruent to $1 \bmod 4$ and occurring with
odd exponent~\cite{e1}. In this article we prove that the number of
supersingular primes is infinite for certain elliptic curves which do
not satisfy any of the above conditions, thereby providing the first
new examples of such curves since the work of Elkies.

Specifically, for $p$ prime, let $w_p$ be the unique Atkin--Lehner
involution~\cite{al} on the modular curve $X_0(p)$, and write
$X_0^*(p)$ for the quotient curve $X_0(p) / w_p$.  Then $X_0^*(p)$ is
a moduli space parameterizing unordered pairs of elliptic curves
$\{E,E'\}$ together with a cyclic $p$--isogeny $\phi\colon E \to E'$.
The main result of this paper is the following:

\begin{theorem}\label{main}
Suppose $p$ is equal to $3$, $5$, $7$, $11$, $13$, or $19$.  Let
$\{E,E'\}$ be a pair of elliptic curves parametrized by a rational
point on the moduli space $X_0^*(p)$, and suppose $E$ does not have
supersingular reduction mod $p$. Then $E$ has infinitely many
supersingular primes.
\end{theorem}

For pairs $E,E'$ whose $j$--invariants are imaginary quadratic
conjugates, the theorem provides new examples of ordinary elliptic
curves with infinitely many supersingular primes. In
Section~\ref{classpoly} we introduce the Heegner point analogues of
Hilbert class polynomials that enable the proof of
Theorem~\ref{main}. Section~\ref{realroots} analyzes the real roots of
these polynomials, and Section~\ref{proof} gives the proof of the
theorem. Section~\ref{numerical} explains the precise relationship
between the curves $E$ of Theorem~\ref{main} and the curves
of~\cite{e1} and~\cite{e2}.

\section{Class polynomial calculations}
\label{classpoly}

Fix an odd prime $p$ such that $X_0^*(p)$ has genus $0$. In this
section we do not impose any other conditions on $p$. Therefore $p$ is
one of $3$, $5$, $7$, $11$, $13$, $17$, $19$, $23$, $29$, $31$, $41$,
$47$, $59$, or $71$. Under these conditions, we will construct a
sequence of polynomials for $X_0^*(p)$ which are analogues to the
Hilbert class polynomials for $X(1)$.  Instead of using CM points on
$X(1)$ we will use Heegner points on $X_0^*(p)$.  We then describe how
our variant class polynomials factor into near perfect squares modulo
primes $\ell \neq p$ and later modulo $\ell = p$.  We also classify
all of the real roots of these polynomials.  Taken together, these
properties of the class polynomials can be used to construct
supersingular primes for points on $X_0^*(p)$.  For $\ell = p$, our
square factorization results only hold for small values of $p$, which
explains why Theorem~\ref{main} is restricted to these values.

The case $p=2$ is omitted because its Heegner points exhibit very
different behavior from the odd case. A discussion of this case can be
found in~\cite{dj}.

For negative integers $D \equiv 0$ or $1 \bmod 4$, write
$\O_D$ for the unique imaginary quadratic order of discriminant
$D$. We assume throughout this chapter that $D$ is of the form
$-p\ell$ or $-4p\ell$ for some prime~$\ell \neq p$. For either choice
of $D$, we denote by $\p$ the ideal of $\O_D$ generated by $p$ and
$\sqrt{D}$.

\begin{lemma}\label{canonical-torsion}
Let $E$ be an elliptic curve over $\C$ with complex multiplication by
$\O_D$. There is exactly one $p$--torsion subgroup of $E$ which is
annihilated by the ideal $\p \subset \O_D$.
\end{lemma}

\begin{proof}
An elliptic curve $E$ with CM by $\O_D$ corresponds to a quotient of
the complex plane $\C$ by a lattice $L$ which is homothetic to an
ideal class in $\O_D$. By scaling $L$ appropriately, we may assume $L
= \<1,\tau\>$ where $\tau = \frac{-b + \sqrt{b^2 - 4ac}}{2a}$ is in
the upper half plane $\H$, with $b^2 - 4ac = D$.

The $p$--torsion subgroups of $E$ are generated in $\C/L$ by
$1/p$, $\tau/p$, $(\tau+1)/p$, \dots, $(\tau+(p-1))/p$. For $z$ to be
annihilated by $\p$ means exactly that the $\R$--linear combination
$\sqrt{D} z = z_1 \cdot 1 + z_2 \cdot \tau$ has integer
coefficients. We have the equations:
\begin{eqnarray}
\label{1/N-torsion-pt}
\sqrt{D} \cdot \frac{1}{p} & = & \frac{b}{p} + \frac{2a}{p} \tau\\
\label{k/N-torsion-pt}
\sqrt{D} \cdot \frac{\tau + k}{p} & = & \frac{bk - 2c}{p} + \frac{2ak
- b}{p} \tau\ \ \ (k = 0,1,2, \ldots, p-1)
\end{eqnarray}

Suppose first that $p \div a$. Then the equation $D = b^2 - 4ac$ means
that $p \div b$, so Equation~\eqref{1/N-torsion-pt} shows that $1/p$
is annihilated by $\p$. By Equation~\eqref{k/N-torsion-pt}, in order
for $(\tau+k)/p$ to be annihilated by $\p$ it would have to be the
case that $p \div (bk - 2c)$, but this cannot happen since $p \div b$
and $p \ndiv c$.

Conversely, if $p \ndiv a$ then Equation~\eqref{1/N-torsion-pt} shows
that $1/p$ is not annihilated by $\p$, and one easily checks using
Equation~\eqref{k/N-torsion-pt} that $(\tau+k)/p$ is annihilated if
and only if $k \equiv b/2a \pmod{p}$.

\end{proof}

One consequence of Lemma~\ref{canonical-torsion} is that, if
$\phi\colon E \to E'$ is the unique cyclic $p$--isogeny whose kernel
is the $p$--torsion subgroup of Lemma~\ref{canonical-torsion}, then
$E'$ also has CM by $\O_D$. Indeed, the lattice $L'$ generated by $L$
and this $p$--torsion subgroup is closed under multiplication by both
$1$ and $\p$, which additively generate all of $\O_D$. In the case
where $D = -p\ell$ and hence $\O_D$ is a maximal order, it follows
immediately that $L'$ has complex multiplication by $\O_D$. When $D =
-4p\ell$, we have to make sure that the CM ring is not an order
strictly containing $\O_{-4p\ell}$, of which the only one is
$\O_{-p\ell}$. But the discriminants of the endomorphism rings of two
$p$--isogenous CM elliptic curves can only differ by a multiple of $p$
if they differ at all~\cite{kohel}, and we have assumed that $p$ is
odd, so the discriminants cannot differ by factors of $2$.

A point on $X_0(p)$ that parameterizes isogenous curves of the same CM
order is called a Heegner point~\cite{gr1}.  We have just showed that
every $E$ with CM by $\O_D$ lifts to a unique Heegner point on
$X_0(p)$.

\begin{definition}\label{heegner}
For any elliptic curve $E$ with CM by $\O_D$, let $\tilde{E}$ denote the
Heegner point on $X_0(p)$ corresponding to the isogeny $E \to E'$
whose kernel is the $p$--torsion subgroup of
Lemma~\ref{canonical-torsion}.
\end{definition}

Let $j_p$ denote a Hauptmodul on $X_0^*(p)$, i.e., a rational
coordinate function on $X_0^*(p)$ with a simple pole of residue $1$ at
$\infty$. Such a function exists since the curve $X_0^*(p)$ always has
a rational cusp and we are assuming its genus is $0$.

\begin{proposition}\label{hN-polynomial}
For each ideal $\a$ of $\O_D$, let $E_\a$ denote the elliptic curve
corresponding to $\C/\a$. For $|D|$ sufficiently large, the minimal
polynomial of $j_p(\tilde{E}_\p)$ over $\Q$ is given by
$$
P_D(X) := \left( \prod_{[\a] \in \Cl(\O_D)} (X - j_p(\tilde{E}_\a))
\right)^{1/2}
$$
where the product is taken over all ideal classes of $\O_D$.
\end{proposition}

\begin{proof}
First, note that $(X - j_p(\tilde{E}_\p))$ is one of the factors in
the product. To get the other factors, start from the known formula
$$
H_D(X) := \prod_{[\a] \in \Cl(\O_D)} (X - j(E_\a))
$$ for the Hilbert class polynomial $H_D(X)$, which by~\cite{cox} is
the minimal polynomial of the $j$--invariant of $E_\p$. Let $G$ be the
absolute Galois group of $\Q$. For every $\sigma \in G$, we have
$\sigma(j(E_\p)) = j(E_{\a})$ for some ideal class $\a$ of $\O_D$
appearing in the above product. We claim that
$\sigma(j_p(\tilde{E}_\p)) = j_p(\tilde{E}_{\a})$ as well, or
equivalently, the map $\sigma\colon E_\p \to E_{\a}$ sends the
distinguished $p$--torsion subgroup of $E_\p$ from
Lemma~\ref{canonical-torsion} to that of $E_{\a}$. But $\sigma$ sends
the endomorphism ring of $E_\p$ into the endomorphism ring of
$E_{\a}$, and in both cases there are only two conjugate embeddings of
$\O_D$ into the endomorphism ring of the elliptic curve, with either
choice resulting in the same action of $\p$ and hence in the same
distinguished $p$--torsion subgroup.

From this claim we see that the set of Galois conjugates of
$j_p(\tilde{E}_\p)$ is exactly $\{j_p(\tilde{E}_\a) \mid \a \subset
\O_D\}$, and so the minimal polynomial contains all the factors in the
product.

We now prove that each linear factor in the product occurs with
multiplicity two. For any ideal class $[\a]$, the Atkin--Lehner image
of $\tilde{E}_\a$ is $\tilde{E}_{\a'}$ for some other ideal class
$[\a'] \in \Cl(\O_D)$ (by the remarks following
Lemma~\ref{canonical-torsion}). The ideal classes $[\a]$ and $[\a']$
are not identical since the $2$-$1$ covering map $\pi\colon X_0(p) \to
X_0^*(p)$ has only finitely many branch points, and we can avoid these
branch points by choosing $|D|$ sufficiently large. Hence
$j_p(\tilde{E}_\a) = j_p(\tilde{E}_{\a'})$, and since $\pi$ is 2 to 1,
these are the only equalities among the roots of the factors in the
product.
\end{proof}

From now on, we will assume that $|D|$ is large enough to
satisfy Proposition~\ref{hN-polynomial}.

\begin{lemma}\label{canonical-reduced}
Let $\P$ be a prime of the splitting field $K$ of $P_D(X)$ lying over
$\ell$, with residue field $k$. Let $\E$ be an elliptic curve
defined over $k$, and fix an embedding
$\O_D \hookrightarrow \End(\E)$. Then there is exactly one $p$--torsion
subgroup of $\E$ which is annihilated by $\p \subset \O_D$.
\end{lemma}
\begin{proof}
By Deuring's Lifting Lemma \cite{deuring}, there is
exactly one lifting of $\E$ to an elliptic curve $E$ over $K$ with CM by $\O_D$
such that reduction mod $\P$ induces the embedding $\O_D
\hookrightarrow \End(\E)$. The $p$--torsion lattices of $E$ and $\E$ are
isomorphic via reduction~\cite{aec}, so the unique $p$--torsion
subgroup of $E$ from Lemma~\ref{canonical-torsion} descends to a
unique $p$--torsion subgroup on $\E$.
\end{proof}

As in Definition~\ref{heegner}, we denote by $\tilde{\E}$ the point on
$X_0(p) \bmod \P$ corresponding to the elliptic curve $\E$ 
together with the cyclic $p$--isogeny whose kernel is the subgroup
determined by Lemma~\ref{canonical-reduced}.

\begin{proposition}\label{hN-mod-l-prelim}
Suppose the odd prime $\ell$ splits in $\O_{-p}$ and $\O_{-4p}$
(equivalently, $-p$ is a quadratic residue modulo $\ell$). Then,
modulo~$\ell$, all roots of the polynomial $P_D(X)$ occur with even
multiplicity, except possibly those corresponding to elliptic curves
with $j \equiv 1728 \bmod \ell$ when $D = -p\ell$, or elliptic curves
which are 2--isogenous to those curves when $D = -4p\ell$.
\end{proposition}
\begin{proof}
Assume first that $D$ is a fundamental discriminant. We show that the
points $\tilde{E}$ corresponding to roots away from $j(E) = 1728$
occur naturally in pairs modulo~$\ell$.  We begin with the following
facts from~\cite{e2} concerning the Hilbert class polynomial $H_D(X)$
defined in the proof of Proposition~\ref{hN-polynomial}.  Each root of
$H_D(X)$ corresponds to an isomorphism class of elliptic curves $E$
with complex multiplication by $\O_D$. The reduction of this root
modulo~$\ell$ corresponds to a reduction of $E$ to a supersingular
elliptic curve $\E$ in characteristic $\ell$, or equivalently an
embedding $\iota\colon \O_D \hookrightarrow
\End(\E)$. Since $\ell$ ramifies in $\O_D$, the conjugate
$\bar{\iota}$ of $\iota$ is again an embedding of $\O_D$ into
$\End(\E)$, and $\E$ lifts by way of $\bar{\iota}$ to an elliptic
curve $E'$ in characteristic zero, which is not isomorphic to $E$
provided that $j(E) \not \equiv 1728 \pmod{\ell}$.

In order to show that the root $j_p(\tilde{\E})$ occurs twice in
$P_D(X)$ modulo~$\ell$, we must show that the two curves $E$ and $E'$
from the previous paragraph correspond to two different roots of
$P_D(X)$ in characteristic zero, and that they both reduce to
$\tilde{\E}$ modulo $l$.  To prove the second fact, observe that the
embeddings $\iota$ and $\bar{\iota}$ both determine the same
$p$--torsion subgroup of $\E$ under Lemma~\ref{canonical-reduced},
since $\p$ equals itself under conjugation, so $\tilde{E}$ and
$\tilde{E}'$ both reduce to $\tilde{\E}$. As for the first fact, we
have $E \neq E'$ provided that $j(E) \not \equiv 1728 \bmod \ell$, so
$\tilde{E} \neq \tilde{E}'$. The only other way $E$ and $E'$ could be
equal on $X_0^*(p)$ is if $w_p(\tilde{E}) = \tilde{E}'$. But if these
two were equal, then in particular their reductions mod~$\ell$ would
be equal, so $w_p(\tilde{\E}) = \tilde{\E}'$. On the other hand, we
have just showed that $\tilde{\E} = \tilde{\E}'$. Putting the two
equations together yields $w_p(\tilde{\E}) = \tilde{\E}$. We show that
this cannot happen.

Let $\phi\colon \E \to \E'$ be the cyclic $p$--isogeny corresponding to
$\tilde{\E}$. The equation $w_p(\tilde{\E}) = \tilde{\E}$ implies that
the dual isogeny $\hat{\phi}$ of $\phi$ is isomorphic to $\phi$, or
that there exist 
isomorphisms $\psi_1\colon \E \to \E'$ and $\psi_2\colon \E' \to \E$ making
the diagram
$$
\xymatrix{
\E \ar[r]^{\phi} \ar[d]_{\psi_1} & \E' \ar[d]^{\psi_2} \\
\E' \ar[r]^{\hat{\phi}} & \E
}
$$ commute. Since
$p$ is prime, the equation 
$\hat{\phi} \phi = [p]$ at once implies that $\psi_2 \phi$ is not equal
to multiplication by any integer, which in turn means that $\psi_2 \phi$
algebraically generates an imaginary quadratic order $\O$ inside
$\End(\E)$. But we also have $(\psi_2 \phi)^2 = u[p]$ for some $u
\in \Aut(\E)$ (specifically, $u = \psi_2 \psi_1$), from which we
conclude that $\O$ contains a square root of $-p$, and thus that $\E$
has CM by either $\O_{-p}$ or $\O_{-4p}$. Moreover, since $\ell$ splits
in these orders by hypothesis, the curve $\E$ must have ordinary
reduction mod $\ell$. On the other hand, by~\cite{deuring} every root
of $H_D(X)$ mod $\ell$
(and hence every root of $P_D(X)$ mod $\ell$)
corresponds to an elliptic curve of supersingular
reduction mod $\ell$, which provides our contradiction.

For the non--fundamental discriminant $D = -4p\ell$, set $D' := -p\ell$ for
convenience. Let $\varepsilon$ be 0, 1, or 2 according as 2 is inert,
ramified or split in~$\O_{D'}$. Then the divisor of zeros $(P_D)_0$ of
$P_D$ in characteristic $\ell$ or $0$ is equal to the Hecke
correspondence $T_2$ on $X_0^*(p)$ applied to the divisor of zeros
$(P_{D'})_0$ of $P_{D'}$, minus $\varepsilon$ times the divisor
$(P_{D'})_0$. That is,
\begin{equation}\label{t2-expansion}
(P_D)_0 = T_2((P_{D'})_0) - \varepsilon (P_{D'})_0.
\end{equation}
Every zero of $P_{D'}$, except for the divisors with $j$--values of
1728, appears in $(P_{D'})_0$ with even coefficient in
characteristic $l$, and hence also appears in $(P_D)_0$ with even
coefficient by~\eqref{t2-expansion}. The only divisors
unaccounted for are those with $j$--values of 1728, and the images of
such divisors under $T_2$, so the Proposition is proved.
\end{proof}

\section{Real roots of $P_D(X)$}\label{real-roots}
\label{realroots}

We find the real roots of the class polynomial $P_D(X)$. 
A real root of $P_D(X)$ corresponds to an unordered pair $\{E, E'\}$
of cyclic $p$--isogenous elliptic curves which is fixed under complex
conjugation. Choose an ideal class $[\a] \in \Cl(\O_D)$ representing
$E$; then $[\a \p]$ represents $E'$. For $\{E,E'\}$ to be fixed
under complex conjugation means that
$$
\{[\a], [\a\p]\} = \{[\conj{\a}], [\conj{\a\p}]\},
$$ where the bar denotes complex conjugation. This can happen in two
ways: either $[\a] = [\conj{\a}]$, or $[\a\p] = [\conj{\a}]$.
\begin{definition}\label{bounded-unbounded}
With notation as above, a real root of $P_D(X)$ is said to be {\em
unbounded} if $[\a] = [\conj{\a}]$, and {\em bounded} if $[\a\p] =
[\conj{\a}]$.
\end{definition}

For the primes $p \equiv 1 \bmod 4$, the behavior of the real roots of
$P_D(X)$ closely mimics the case of $H_D(X)$ which was treated
in~\cite{e1}. This is not surprising given that $X(1) = X_0(1)$ is a
special case of $X_0(p)$ when $p \equiv 1 \bmod 4$. However, when $p
\equiv 3 \bmod 4$ the real roots of $P_D(X)$ exhibit very different
behavior. It is therefore necessary to treat the two cases separately.

\subsection{The case $p \equiv 1 \bmod 4$}
\label{real-roots-1mod4}

In this section, we assume that $p \equiv 1 \pmod{4}$ and that $D$ is
equal to $-p\ell$ or $-4p\ell$, where $\ell$ is chosen to be a prime
congruent to $3 \bmod 4$ which splits in $\O_{-p}$ and $\O_{-4p}$.

An unbounded real root of $P_D(X)$ corresponds to an isogeny $E \to
E'$ which is isomorphic to itself under complex conjugation, meaning
that $\tilde{E}$ is a real point on $X_0(p)$.  Since the covering
$X_0(p) \to X(1)$ is defined over $\Q$, each such real point $\tilde{E}$
has $j(E)$ real, so we can count these points by counting the ideal
classes $[\a]$ for which $j(\a)$ is real. By genus theory~\cite{cox},
there are two such ideal classes for $\O_{-p l}$ and two for $\O_{-4p
l}$, corresponding to the quadratic forms
\begin{eqnarray*}
& x^2 + xy + \left(\frac{p\ell+1}{4}\right)y^2 & \\
& px^2 + pxy + \left(\frac{p+\ell}{4}\right)y^2 &
\end{eqnarray*}
for $D = -p\ell$, and
\begin{eqnarray*}
& x^2 + p\ell y^2 & \\
& p x^2 + \ell y^2 &
\end{eqnarray*}
for $D = -4p\ell$.

Since the first two forms above are Atkin--Lehner images of each
other, and the last two are Atkin--Lehner images of each other, the
first pair of real points on $X_0(p)$, upon quotienting by $w_p$,
yields one real root of $P_{-p\ell}(X)$, and the second pair yields a
real root for $P_{-4p\ell}(X)$. For $D = -p\ell$, the quadratic form
$x^2 + xy + \left(\frac{p\ell+1}{4}\right)y^2$ has the root 
$\tau = (-1 + \sqrt{-p\ell})/2$ in the upper half plane, and
$$
\lim_{\ell \to \infty} j_p\left(\frac{-1+\sqrt{-p\ell}}{2}\right) =
-\infty.
$$
Similarly, for $D = -4 p\ell$,
the quadratic form $x^2 + p\ell y^2$ has the root $\tau =
\sqrt{-p\ell}$ with $\lim_{\ell \to \infty} j_p(\tau) = \infty$. The
divergence of the roots $j_p(\tau)$ of $P_D(X)$, as $\ell \to \infty$,
justifies the terminology ``unbounded.''

A bounded real root of $P_D(X)$ occurs when $[\a\p] = [\conj{\a}]$, or
equivalently $[\p] = [\conj{\a}]^2$. Viewing each ideal class as a
quadratic form, a bounded root exists if and only if the quadratic
form $px^2 + \ell y^2$ (for $D=-4p\ell$) or $px^2 + pxy +
\frac{p+\ell}{4} y^2$ (for $D=-p\ell$) is equal to the direct
composition of some quadratic form $ax^2 + bxy + cy^2$ with itself. In
particular, this implies by definition of composition that there
exists a nonzero integer $z$ satisfying the Diophantine equation $px^2
+ \ell y^2 = z^2$ in the $D = -4p\ell$ case, or $px^2 + pxy +
\frac{p+\ell}{4}y^2 = z^2$ in the $D = -p\ell$ case. We show that this
cannot happen in our situation.
\begin{lemma}
The Diophantine equations $px^2 + \ell y^2 = z^2$ and $px^2 + pxy +
\frac{p+\ell}{4}y^2 = z^2$ have no nonzero solutions $x,y,z \in \Z$.
\end{lemma}
\begin{proof}
Suppose there were a nonzero solution. We may assume $y \not \equiv 0
\pmod{p}$, or else descent yields a contradiction. Then reducing the
equations modulo $p$, we get that $\ell$ is a perfect square mod $p$,
which contradicts the assumptions that $\ell \equiv 3 \bmod 4$ and $l$
splits in $\O_{-p}$.
\end{proof}

We conclude that the polynomial $P_D(X)$ has one unbounded real root
and no bounded real roots, with the bounded real root diverging to
$\infty$ for $D = -4p\ell$ and $-\infty$ for $D = -p\ell$, as $\ell
\to \infty$.

\subsection{The case $p \equiv 3 \bmod 4$}
\label{real-roots-3mod4}

We assume that $p \equiv 3 \pmod{4}$ and that $D = -4p\ell$, where
$\ell \equiv 1 \bmod 4$ and $\ell$ splits in $\O_{-p}$ and $\O_{-4p}$.
Using genus theory as before, the unbounded real root of $P_D(X)$ is
represented by the pair of ideals corresponding to the two
$w_p$--equivalent quadratic forms
\begin{eqnarray*}
& x^2 + p\ell y^2 & \\
& p x^2 + \ell y^2. &
\end{eqnarray*}
Hence the polynomial $P_D(X)$ has one unbounded real root,
which approaches $\infty$ as $\ell$ becomes large.

Recall that a bounded real root corresponds to an equivalence class of
quadratic forms $ax^2+bxy+cy^2$ whose square in the form class group
is equal to the form $px^2 + \ell y^2$.  There is at most one such
form class, because a second one would result in more $2$--torsion
classes in the ideal class group of $\O_D$ than were found in the
preceding analysis of the unbounded roots.

To show the existence of such a quadratic form, it suffices to
construct a quadratic form $pax^2+bxy+ay^2$ of discriminant $D$ with
$p$ dividing $b$. Indeed, the Dirichlet composition~\cite{cox} of
$pax^2 + bxy + ay^2$ with itself is $a^2x^2 + bxy + py^2$, which is
properly equivalent to $px^2 + \ell y^2$ since $p \div b$ and the
discriminants of the two forms match.

To find such a quadratic form, choose integers $A$ and $B$ such that
$\ell = A^2 - p B^2 = (A + B \sqrt{p}) (A - B \sqrt{p})$. Such
integers exist because $\ell$ splits in $\Q(\sqrt{p})$, and all such
representations of $\ell$ differ by a factor of $\pm \varepsilon^n$
where $\varepsilon := c + d \sqrt{p}$ is the fundamental unit of
$\Q(\sqrt{p})$. Note that $c$ and $d$ are integers, since $p \equiv 3
\bmod 4$, and that $c$ is even and $d$ is odd. Accordingly,
multiplication by $\varepsilon$ changes the parity of $A$, so there
exist representations with $A$ even and with $A$ odd. Choose $A$ to be
odd, and set $a=A$, $b=2pB$ to obtain a quadratic form $pax^2 + bxy +
ay^2$ of discriminant $-4p\ell$.

We now find the minimal possible value for $B$ (equivalently, the
minimal possible $b$), subject to the constraint that $A$ is odd. This
value for $B$ is determined by the requirement that multiplication by
$\varepsilon^2$ must increase the size of the coefficients of the
factor $A - B \sqrt{p}$. We compute these coefficients to be:
$$
(A - B \sqrt{p}) (c + d \sqrt{p})^2 = (A c^2 - 2 B cd p + A d^2 p) +
(2Acd - Bc^2 - Bd^2 p) \sqrt{p}
$$

The requirement is thus $B < (2Acd - Bc^2 - Bd^2 p)$, or
$$
\frac{B}{A} < \frac{2 cd} {c^2 + d^2 p + 1} = \frac{d}{c} \cdot
\frac{2c^2}{c^2 + d^2p + 1}.
$$
But $d^2p = c^2 - 1$, so the fraction $(2c^2)/(c^2+d^2p+1)$ equals
$1$, whence our condition on $B$ is just $B/A < d/c$. One could
have done the same computation using the inequality on $A$ given by
the other coefficient; the reader can verify that doing so produces
the same inequality.

Now, if $b$ is chosen to be minimal and of the above form (i.e., $cB <
dA$, or equivalently $cb < 2pda$, and $A$ is odd), then the root $\tau
= \frac{-b+\sqrt{D}}{2pa}$ of the quadratic form $pax^2 + bxy + ay^2$
lying in the upper half plane has absolute value $1/\sqrt{p}$ and real
part equal to $-B/A$, with $-d/c < -B/A < 0$. Denote the set of all
such complex numbers in the upper half plane by $S$.  Since all of the
points on the circular arc $S$ are distinct in $X_0^*(p)$, the
function $j_p(z)$ is monotonic (and, of course, real valued) in the
clockwise direction along this circular arc. From $q$--expansions we
see that $j_p(z)$ is in fact increasing clockwise along the arc
$S$. We claim that, for random large values of $\ell$, the locations
of the corresponding roots $\tau$ (as a function of $\ell$) are
uniformly distributed along the arc $S$ in a weak sense to be made
precise in Lemma~\ref{uniform-dist}. It follows that the bounded
real root of the polynomial $P_D(X)$ is uniformly distributed along
the real interval $j_p(S)$ as $D$ varies.

\begin{lemma}\label{uniform-dist}
Let $\A$ be an arithmetic progression containing infinitely many
primes $\ell$ which are congruent to $3 \bmod 4$ and split in
$\O_{-p}$ and $\O_{-4p}$. For any sub-arc \mbox{$T \subset S$} of
nonzero length, there exist infinitely many primes $\ell \in \A$ whose
corresponding roots $\tau$ above lie in $T$.
\end{lemma}
\begin{proof}
Let $U$ be the projection of $T$ to the real axis. Using the fact that
$\Re(\tau) = -B/A$, we see that it suffices to show that $-B/A \in U$
for infinitely many primes $\ell \in \A$.  Consider the function
$$
\sigma(\a) := \left(\frac{N(\a)}{N(\conj{\a})}\right)^{\frac{2
    \pi i}{\log(\varepsilon/\conj{\varepsilon})}}
$$ mapping ideals $\a$ of $\O_p$ into complex numbers of norm $1$. Let
$\sigma(A,B)$ denote the value of $\sigma$ on the principal ideal
$(A+B\sqrt{p})$ in $\O_p$. Then $\sigma(A,B)$ is purely a function of
$B/A$, and as $B/A$ increases from $0$ to $d/c$ with $B$ positive, the
point $\sigma(A,B) \in S^1$ increases monotonically in angle from $0$
to $2\pi$. Thus it is enough to show that $\sigma(A,B)$ is
equidistributed on $S^1$ where $A,B$ vary as a function of $\ell \in
\A$, with $\ell = A^2 - pB^2$. But the equidistribution of values of
$\sigma$ with respect to $\ell$ has already been proven
in~\cite[p. 318]{lang}.
\end{proof}

In summary, for $p \equiv 3 \bmod 4$ and $D = -4p\ell$, where $\ell
\equiv 1 \bmod 4$ and $\ell$ splits in $\O_{-p}$ and $\O_{-4p}$, the
polynomial $P_D(X)$ has exactly two real roots, with the unbounded
real root diverging to $\infty$ as $\ell$ increases and the bounded
real root being uniformly distributed in the real interval $j_p(S)$ as
the prime $\ell$ is varied.

\section{Proof of the main theorem}
\label{proof}

\subsection{Specification of Hauptmoduls}
\label{concrete}

For the sake of concreteness, we will use the following Hauptmoduls
for the curves $X_0^*(p)$, $p = 3,5,7,11,13,19$. The derivation of
these Hauptmoduls is discussed in~\cite{e3}.

For $p=3,5,7,13$, the modular curve $X_0(p)$ is a rational curve with
coordinate
\begin{equation}\label{j30}
j_{p,0}(z) := \left( \frac{\eta(z)}{\eta(pz)} \right)^{\frac{24}{p-1}},
\end{equation}
where $\eta$ is the Dedekind eta function. The action of the
Atkin--Lehner involution $w_p$ is given by

\begin{equation}\label{wp-jp}
w_p(j_{p,0}(z)) = \frac{p^{\frac{12}{p-1}}}{j_{p,0}(z)}.
\end{equation}
For these primes, we use the Hauptmodul $j_p$ defined by the formula
\begin{equation}\label{jp-def}
j_p(z) := j_{p,0}(z) + w_p(j_{p,0}(z)).
\end{equation}

For $p=11$ we use the Hauptmodul
$$
j_{11}(z) := \left( \frac{\theta_{1,1,3}(z)}{\eta(z) \eta(11 z)}
\right)^2,
$$
where $\theta_{a,b,c}(z)$ is defined to be the theta function
$$
\theta_{a,b,c}(z) := \sum_{x,y \in \Z} q^{ax^2 + bxy + cy^2},\ \ \
q:=e^{2 \pi i z},
$$
valid for all $z$ in the upper half plane $\H$.

For $p=19$, we use the function
$$
j_{19}(z) := \left( \frac{\theta_{1,1,5}(z)}{\theta^*_{1,1,5}(z)}
\right)^2,
$$
where now $\theta^*_{1,1,5}(z)$ is defined by
$$
\theta^*_{1,1,5}(z) := \sum_{m+n\equiv 1 (2)} (-1)^m q^{\frac{1}{2}
(m^2 + mn + 5n^2)},\ \ \ q:=e^{2\pi i z}.
$$

\subsection{Proof of the theorem for $p \equiv 3 \bmod 4$}
\label{3mod4main}
We assume that $p$ is equal to $3$, $7$, $11$, or $19$.  As before, we
will use the polynomials $P_D(X)$, $D = -p\ell$ or $D = -4p\ell$,
where the prime $\ell$ is both $1 \bmod 4$ and a quadratic residue mod
$p$. Note that $P_D(X)$ is monic (since its roots $j_p(\tilde{E})$ are
algebraic integers) and each such curve $E$ is supersingular mod $p$
and mod $\ell$ (since $p$ and $\ell$ ramify in $D$).

\begin{proposition}\label{hN-mod-l}
The polynomial $P_D(X)$ is a square modulo $\ell$.
\end{proposition}
\begin{proof}
By Proposition~\ref{hN-mod-l-prelim}, we only have to exclude the
possibility of there being roots associated to the $j$--invariant
$1728$. First consider the case $D = -p\ell$. Suppose $j_p(\tilde{E})$
were a root of $P_D(X)$, with $j(E) \equiv 1728 \bmod \ell$. Then $E$
would be supersingular mod $\ell$ and have complex multiplication by
$\O_{-4}$. But $\ell$ splits in $\O_{-4}$, so a curve with CM by
$\O_{-4}$ cannot be supersingular mod $\ell$.

Now take $D = -4p\ell$. As in the proof of
Proposition~\ref{hN-mod-l-prelim}, set $D' := -p\ell$. Then, since all
coefficients in the divisor of zeros $(P_{D'})_0$ are even in
characteristic $\ell$, the proof of Proposition~\ref{hN-mod-l-prelim}
shows that every coefficient of $(P_D)_0$ is even as well.
\end{proof}

\begin{lemma}\label{hN-mod-N}
For $D = -4p\ell$, the polynomial $P_D(X)$ is a perfect square modulo
$p$.
\end{lemma}
\begin{proof}
Suppose first that $p = 3$ or $7$. Every root of $P_D(X)$ is of the
form $j_p(\tilde{E})$ where $E$ is a supersingular elliptic curve mod
$p$. But there is only one isomorphism class of supersingular elliptic
curves mod $p$. It follows that $P_D(X)$ has divisor of zeros equal to
$\deg(P_D) \cdot (j_p(\tilde{E}))$. Since $P_D(X) \bmod p$ is monic,
has even degree, and has only one root of maximal multiplicity, it
must be a perfect square.

Now suppose $p = 11$. Write $D' = -p\ell$ as before. Here there are
two isomorphism classes of supersingular elliptic curves mod $p$,
having the values $0$ and $-1$ under the coordinate function $j_{11}$
of \S\ref{concrete}. Using the algorithm of Pizer~\cite{pizer}, we
find that the action of the Hecke correspondence $T_2$, as given by the
Brandt matrix $B(2)$, is represented by the equations:
\begin{eqnarray*}
T_2((0)) & = & 1\cdot(0) + 2\cdot(-1) \\
T_2((-1)) & = & 3\cdot(0) + 0\cdot(-1).
\end{eqnarray*}
Since the roots of the polynomial $P_{D'}(X)$ are supersingular, the
polynomial $P_{D'}(X)$ has the form $X^m (X+1)^n \bmod 11$ for some
integers $m$ and $n$. The above calculation of $T_2$, combined with
Equation~\eqref{t2-expansion}, yields
$$
P_D(X) \equiv X^{m+3n - \varepsilon m} (X+1)^{2m - \varepsilon n}
\bmod 11,
$$ which is a perfect square since $\deg(P_{D'}) = m+n$ is even and
$\varepsilon$ is even for all primes $\ell \equiv 1 \bmod 4$ which are
squares $\bmod p$.

The case $p=19$ is similar: using the Hauptmodul $j_{19}$ of
\S\ref{concrete}, the Hecke correspondence mod $19$ has matrix
$\left[\begin{smallmatrix}1 & 2 \\ 1 & 2\end{smallmatrix}\right]$ with
respect to the basis of supersingular invariants $\{(0), (8)\}$. Since
the columns of this matrix add up to even numbers, the polynomial
$P_D(X)$ is always a perfect square modulo $19$ for $D = -4p\ell$ and
our choices of $\ell$.
\end{proof}

\begin{theorem}\label{main-N}
Suppose $p = 3,7,11,$ or $19$. Let $\{E,E'\}$ be a pair of elliptic
curves, defined over $K$, corresponding to a rational point on
$X_0^*(p)$, and assume that $E$ is not supersingular at $p$. Then $E$
has infinitely many supersingular primes.
\end{theorem}

\begin{proof}
If $E$ is represented by the complex lattice $\<1,\tau\>$ with $\tau
\in \H$, the fact that $h := j_p(\tau)$ is real means that we may
(c.f. Section~\ref{real-roots}) take $\tau$ either on the unbounded
arcs corresponding to $\Re(\tau) = 0$ or $\Re(\tau) = 1/2$, or on the
bounded arc $j_p(S)$ of Lemma~\ref{uniform-dist}.  In the unbounded
case, $j(\tau)$ is real and $K$ has a real embedding, so the result
follows from~\cite{e2} and we do not need to do it here. We can
therefore assume that $\tau \in S$ and $-d/c \leq \Re(\tau) \leq
0$. Moreover, we can assume these inequalities are strict, since
otherwise $E$ has CM and its supersingular primes are known to have
density $1/2$.

Now suppose $h$ is rational inside the interior of the interval
$j_p(S)$ and the curve $E$ is not supersingular modulo $p$. Given any
finite set $\Sigma$ of primes of $K$, we construct a supersingular
prime $\pi$ of $E$ outside of $\Sigma$.

Without loss of generality, suppose that $\Sigma$ contains all of the
primes of bad reduction of $E$. Choose a large prime $\ell$ such that
\begin{enumerate}
\item\label{cond1} $\ell \equiv 1 \bmod 4$ and $\ell$ splits in
  $\O_{-p}$ and $\O_{-4p}$.
\item\label{cond2} $\jacobi{v}{p\ell} = 1$ for every rational prime
  $v$ lying under a prime in $\Sigma$, except possibly $v=p$.
\item\label{cond3} $P_D(h) < 0$.
\end{enumerate}
Condition~\ref{cond3} is satisfied as long as the bounded root $r$ of
$P_D(h)$ falls to the left of $h$ on the real line. Since $h$ is not
on the boundary of $j_p(S)$, Lemma~\ref{uniform-dist} assures the
existence of infinitely many primes $\ell$ satisfying all the
conditions.

The rational number $P_D(h)$ is then congruent to a perfect square
mod $\ell$ (by Proposition~\ref{hN-mod-l}) and mod $p$ (by
Lemma~\ref{hN-mod-N}). However, being negative, it also contains
a factor of $-1$, which is not a perfect square mod $p\ell$. Therefore at
least one of its prime factors $q$ satisfies the equation
\mbox{$\jacobi{q}{p\ell} \neq 1$} and thus is ramified or inert in
$\Q(\sqrt{D})$. Moreover, the denominator of $P_D(h)$ is a perfect
square, since $P_D(X)$ is monic with integer coefficients and
even degree.  Hence we may take $q$ to be a factor of the numerator of
$P_D(h)$. 
Furthermore, $q$ is not equal to $p$, because by hypothesis $E$ is
not supersingular at $p$ so $p$ cannot divide $P_D(h)$.

It follows from Condition~\ref{cond2} that $q$ does not lie under any
prime in $\Sigma$, and $h$ is a root of $P_D(X)$ in characteristic
$q$. Therefore $j(E)$ is a root of $H_D(X)$ in characteristic
$q$. Hence, for any prime $\q$ of $K$ lying over $q$, the
reduction of $E$ at $\q$ has complex multiplication by $\O_{D'}$ for
some factor $D'$ of $D$ such that $D/D'$ is a square, and since $q$ is
not split in $\Q(\sqrt{D})$, it follows that there is a new
supersingular prime $\pi \not \in \Sigma$ lying above $q$.
\end{proof}

\subsection{Proof of the theorem for $p \equiv 1 \bmod 4$}
\label{1mod4}

We now give a proof of Theorem~\ref{main} for the primes $p=5$ and
$13$. Let $\ell$ be a prime congruent to $3 \bmod 4$ such that
$\ell$ splits in $\O_{-p}$ and $\O_{-4p}$. Explicitly, $\ell
\equiv 3,7 \pmod{20}$ for $N=5$, and $\ell \equiv 7,11,15,19,31,47
\pmod{52}$ for $N=13$. Note that Proposition~\ref{hN-mod-l-prelim}
applies in this case. Throughout this section we will use the
Hauptmoduls $j_5$ and $j_{13}$ specified in \S\ref{concrete}.

\begin{proposition}\label{h5-mod-l}
For $p=5$ and $D = -p\ell$ or $D=-4p\ell$, the polynomial $P_D(X)$ is
of the form $(X+22) R(X)^2$ modulo $\ell$.
\end{proposition}
\begin{proof}
From class number considerations we know that the class polynomial
$P_D(X)$ has odd degree. We show that the only factors of $P_D(X)$
lying over $j = 1728$ are equal to $(X+22)$ mod $\ell$. This will imply
that our polynomial has the required form, by
Proposition~\ref{hN-mod-l-prelim}.

Let $E = \C/L$ where $L = \Z[i]$. Then $j(E) = 1728$ and there are six
points (counting multiplicity) of $X_0(5)$ lying over $E$. We compute
the values under $j_{5,0}$ and $j_5$ for each of the choices of
$5$--torsion subgroup of $E$:
\begin{center}
\begin{tabular}{c|c|c}
subgroup & $j_{5,0}$ & $j_5$ \\ \hline
$\<1/5\>$ & $125 + 2 \sqrt{5}$ & $248 + 126 \sqrt{5}$ \\
$\<i/5\>$ & $125 + 2 \sqrt{5}$ & $248 + 126 \sqrt{5}$ \\
$\<(i+1)/5\>$ & $125 - 2 \sqrt{5}$ & $248 - 126 \sqrt{5}$ \\
$\<(i-1)/5\>$ & $125 - 2 \sqrt{5}$ & $248 - 126 \sqrt{5}$ \\
$\<(i+2)/5\>$ & $-11 + 2 i$ & $-22$ \\
$\<(i-2)/5\>$ & $-11 - 2 i$ & $-22$
\end{tabular}
\end{center}
Notice that the two subgroups $G$ of $E$ with $j_5(E,E/G) = -22$ are
characterized by the property $G = iG$. We will use this
characterization to prove that the roots of $P_D(X)$ over $1728$ must
have $j_5 = -22$.

Suppose first that $D = -5\ell$ is a fundamental discriminant. Let
$j_5(\tilde{E})$ be a root of $P_D(X)$ modulo $\ell$ with $j(E) =
1728$ modulo $\ell$. Then the reduction $\E$ of $E$ modulo $\ell$ has
quaternionic endomorphism ring $A$ containing a subring generated by
$\Z[I, (D + \sqrt{D})/2]$, where $I^2 = -1$ and $\sqrt{D}$ in $A$ is
obtained from the embedding $\iota\colon \O_D \to A$ induced by the
reduction map from $E$ to $\E$. Now, the reduction of the ring $A$
modulo 5 is isomorphic to $M_{2 \times 2}(\Z/5)$, with the isomorphism
being given by the action of $A$ on the $5$--torsion group $E[5] =
\E[5]$ of $E$. The element $\sqrt{D}$ has square equal to $D \equiv 0
\bmod 5$, so it is nilpotent in $M_{2\times 2}(\Z/5)$ with kernel
equal to $\ker(5,\sqrt{D}) = \ker \p$. Observe that $\ker(I \sqrt{D}
I^{-1}) = I \ker(\sqrt{D}) = I \ker \p$; on the other hand, $\ker(I
\sqrt{D} I^{-1}) = \ker(\bar{\iota}(\sqrt{D})) = \ker(\sqrt{D}) = \ker
\p$. Therefore the distinguished $5$--torsion subgroup $G = \ker
\p$ of Lemma~\ref{canonical-reduced} satisfies the equality $G
= iG$, as desired.

For the non--fundamental case $D = -20\ell$, note that the Hecke
correspondence $T_2$ applied to the value $j_5(\tilde{E}) = -22$ is a
formal sum of terms all with even coefficient except for $-22$ itself,
so by the proof of Proposition~\ref{hN-mod-l-prelim}, the polynomial
$P_D^5(X)$ is a perfect square except for a linear factor of $(X+22)$.
\end{proof}

\begin{proposition}\label{h13-mod-l}
For $p=13$ and $D = -p\ell$ or $D = -4p\ell$, the polynomial $P_D(X)$
is of the form $(X+6) R(X)^2$ modulo $\ell$.
\end{proposition}
\begin{proof}
Let $j(E) = 1728$. By the same proof as in Proposition~\ref{h5-mod-l},
the kernel $G$ of $\tilde{E}$ satisfies $G = iG$. There are only two
13--torsion subgroups $G$ of $\C/\Z[i]$ that satisfy the equation $G =
iG$, and they are generated respectively by $(2+3i)/13$ and
$(3+2i)/13$. One calculates that $j_{13,0} = -3\pm 2i$ and $j_{13} = -6$
for these points, so as in Proposition~\ref{h5-mod-l} the polynomial
$P_D^{13}$ factors as $(X+6)$ times a perfect square.
\end{proof}

Define $P_\ell(X)$ to be the monic polynomial $P_{-p\ell}(X) \cdot
P_{-4p\ell}(X)$. Then, by Propositions~\ref{h5-mod-l}
and~\ref{h13-mod-l}, the polynomial $P_\ell(X)$ is a perfect square
mod~$\ell$, and by the classification of the real roots of $P_D(X)$ in
\S\ref{real-roots-1mod4}, the polynomial $P_\ell(X)$ has exactly two
real roots which diverge to infinity in opposite directions as $\ell
\to \infty$. In particular, for any fixed real number $h$, the value
of $P_\ell(h)$ is negative for all sufficiently large $\ell$.

\begin{lemma}\label{hN-mod-N-1mod4}
For $p=5$ or $13$, the polynomial $P_\ell(X)$ is a perfect square
modulo $p$.
\end{lemma}
\begin{proof}
Since the polynomial has even degree, it suffices to prove that all
roots of the polynomial are congruent mod~$p$. But every root of
$P_\ell(X) \bmod p$ is of the form $j_p(\tilde{E})$ where $\tilde{E}$
is an elliptic curve whose reduction modulo $p$ is supersingular. For
either $p=5$ or $p=13$, there is only one isomorphism class of
supersingular $j$--invariants mod~$p$, so all such curves $E$ are
isomorphic mod~$p$ and they all have the same $j_p$ value.
\end{proof}

\begin{theorem}\label{main-1mod4}
Suppose $p$ equals $5$ or $13$. Let $\{E,E'\}$ be a pair of elliptic
curves, defined over a number field $K$, corresponding to a rational
point on the curve $X_0^*(p)$. Assume that $E$ is not supersingular at
$p$. Then $E$ has infinitely many supersingular primes.
\end{theorem}
\begin{proof}
Suppose $h := j_p(\tilde{E})$ is rational and not of supersingular
reduction modulo $p$. Given any finite set $\Sigma$ of primes of $K$,
containing all of $E$'s primes of bad reduction, we construct a
supersingular prime $\pi$ of $E$ outside of $\Sigma$.

Choose a large prime $\ell$ satisfying the conditions:
\begin{enumerate}
\item $\ell \equiv 3 \bmod 4$ and $\ell$ splits in $\O_{-p}$ and
  $\O_{-4p}$.
\item $\jacobi{v}{p\ell} = 1$ for every rational prime $v$ lying under
  a prime in $\Sigma$ (except possibly $v=p$).
\item $P_\ell(h) < 0$.
\end{enumerate}
Then the numerator $z$ of the rational number $P_\ell(h)$ is divisible
by some rational prime $q$ which is ramified or inert in
$\Q(\sqrt{D})$ for one of $D = -p\ell$ or $D = -4p\ell$ (equivalently,
has $\jacobi{q}{p\ell} \neq 1$). Indeed, if not, then the absolute
values of both the numerator and the denominator of $P_\ell(h)$ would
have quadratic character 1 modulo $p\ell$. But $\jacobi{-1}{p\ell} =
-1$ by our choice of $\ell$, so the number $P_\ell(h)$ itself would
have quadratic character $-1$ modulo $p\ell$, contradicting the fact
that $P_\ell(h)$ is a perfect square mod $p$ and mod $\ell$.

Moreover, $q$ is not equal to $p$, since the assumption that $E$ is
not supersingular at $p$ implies that $p$ does not divide $P_\ell(h)$.

It follows that $q$ does not lie under any prime in $\Sigma$, and $h$
is a root of $P_\ell(X)$ in characteristic $q$. Therefore, for one of
$D = -p\ell$ or $D = -4p\ell$, the value $j(E)$ is a root of
$H_D(X)$ in characteristic $q$. Hence, for any prime $\q$ of $K$
lying over $q$, the reduction $E_\q$ has complex multiplication by
$\O_{D'}$ for some factor $D'$ of $D$ such that $D/D'$ is a square,
and since $q$ is not split in $\Q(\sqrt{D})$, it follows that there is
a new supersingular prime $\pi \not \in \Sigma$ lying above $q$.
\end{proof}

\section{Numerical Computations}
\label{numerical}
\subsection{Relationship to Elkies's work}

In addition to proving the infinitude of supersingular primes for
elliptic curves defined over real number fields in~\cite{e2}, Elkies
in~\cite[p. 566]{e1} notes that his method also works for
$j$--invariants ``such that the exponent of some prime congruent to
$+1 \bmod 4$ in the absolute norm of $j-12^3$ is odd.''  Thus, even
for the case of elliptic curves over imaginary number fields our
results do not represent the first demonstration of infinitely many
supersingular primes for ordinary curves.  However, one can prove by
direct computation that, over non--real number fields, the set of
elliptic curves given in the statement of Theorem~\ref{main} is
disjoint from the set of curves which satisfy the property stated by
Elkies above.  As an illustration of this fact we will perform the
computation for the case of $X_0^*(3)$.

We preserve the notation from Section~\ref{concrete}. We will need the
equation
\begin{equation}\label{j3-to-j}
j(z) - 1728 = \frac{(j_{3,0}(z)^2 - 486\,j_{3,0}(z) -
  19683)^2}{j_{3,0}(z)^3},
\end{equation}
obtained as in~\cite{e3} by linear algebra on the Fourier coefficients
of $q$--expansions. Because~\cite{e2} already provides for the case of
elliptic curves with real $j$--invariants, we are interested only in
the case of non--real $j$--invariants.  Equations~\eqref{jp-def}
and~\eqref{j3-to-j} show that the only way a rational number $j_3(z)$
can arise from a non--real number $j(z)$ is if the two complex numbers
$j_{3,0}(z)$ and $w_3(j_{3,0}(z))$ are imaginary quadratic complex
conjugates of each other.  When this happens, Equation~\eqref{wp-jp}
then shows that the two complex conjugates multiply to $3^6$, so we
conclude that the norm of $j_{3,0}(z)$ must equal $3^6$.

Taking the norms of both sides of~\eqref{j3-to-j}, we get
\begin{eqnarray*}
\No(j(z)-1728) & = & \frac{\No(j_{3,0}(z)^2 - 486\,j_{3,0}(z) -
  19683)^2}{\No(j_{3,0}(z))^3} \\
& = & \frac{\No(j_{3,0}(z)^2 - 486\,j_{3,0}(z) -
  19683)^2}{(3^6)^3},
\end{eqnarray*}
where the last equality follows from the fact that $j_{3,0}(z)$ has
norm $3^6$. This equation shows that the rational number
$\No(j(z)-1728)$ is always a perfect square, and hence it cannot
satisfy the requirement of Elkies that it possess a prime factor of
odd multiplicity.

\subsection{Points on $X_0^*(11)$}

For a numerical demonstration of our supersingular prime finding
algorithm, consider the point $j_{11} = \frac{21}{2}$ on $X_0^*(11)$,
having $j$--invariant
$$j = \frac{-489229980611 - 42355313 \sqrt{-84567}}{4096},$$ with
$\No(j-1728) = (7646751287/64)^2$.  We find supersingular primes for
this $j$--invariant using class polynomials on $X_0^*(11)$.  For this
we must pick primes $\ell \equiv 1 \bmod 4$ such that $\ell$ is a
quadratic residue mod $11$ and the class polynomial of discriminant
$-44\ell$ has a real root to the left of $\frac{21}{2}$, in order to
ensure that $P_D(\frac{21}{2})$ is negative.

Using $\ell=5$, we find that 
$$
P_{-220}(X) = X^2 - 77 X + 121.
$$
The rational number $P_{-220}(\frac{21}{2}) = -2309/4$ is
negative and a perfect square modulo $55$, so the prime factor $2309$
in the numerator is a supersingular prime for this point.

To find a new supersingular prime not equal to $2309$, we need a new
value of $\ell$ with $\jacobi{2309}{11\ell} = 1$. Using $\ell=37$, we
have
\begin{eqnarray*}
P_{-1628}(X) & = & X^8 - 101042 X^7 - 2728753 X^6 - 167281605 X^5
+ \\
& & 1453552981 X^4 - 4464256335 X^3 + 8630555868 X^2 - \\
& & 9354295951 X + 4253517961
\end{eqnarray*}
and
$$
P_{-1628}({\textstyle \frac{21}{2}}) = -\frac{7^2 \cdot 151 \cdot
  452233314041}{256}.
$$
Of the primes in the numerator, both $7$ and $151$ are quadratic
nonresidues mod $11\cdot 37 = 407$, so our $j$--invariant is
supersingular modulo these primes.  In this case the primes are small
enough to check directly against the tables of supersingular
$j$--invariants in~\cite{ant}; thus we find that $(-489229980611 -
42355313 \sqrt{-84567})/4096$ is congruent to $6 \bmod 7$, and to $67
\bmod 151$ (or to $101 \bmod 151$ if the other square root is chosen),
and that these values are indeed supersingular invariants modulo $7$
and $151$ respectively.

\section{Further directions}
\label{further}

The proofs given here are not limited to the case where $j_p(E)$ is
rational. When $p \equiv 1 \bmod 4$, we can generalize
Theorem~\ref{main} to the case of elliptic curves $E$ whose
$j_p$--invariant has odd algebraic degree. The proof is the same as
that given in~\cite{e1}: for large enough values of $\ell$, the
absolute norm of $P_\ell(j_p(E))$ is negative and hence has a prime
factor lifting to a new supersingular prime of $E$. Likewise, for $p
\equiv 3 \bmod 4$, we can extend our proof to all curves $E$ for which
$j_p(E)$ is real. In this case we assume that all the real conjugates
of $j_p(E)$ lie inside the set $j_p(S)$ of Lemma~\ref{uniform-dist},
since otherwise we can use~\cite{e2} directly. Because the bounded
root of $P_D(X)$ is uniformly distributed along $j_p(S)$, there exists
a value of $D$ making $P_D(X)$ negative valued on exactly one real
conjugate of $j_p(E)$. For this choice of $D$, the numerator of
the absolute norm of $P_D(j_p(E))$ produces a new supersingular prime
for $E$.

One might naturally ask how to prove Theorem~\ref{main} for the primes
$p=17$ or $p>19$. Our proof relies on the fact that the polynomial
$P_D(X)$ is a square mod $p$. When $X_0(p)$ has genus $0$, this fact
is automatic since $P_D(X)$ has only one root in characteristic $p$.
For the genus $1$ cases $p=11$ and $19$, we proved squareness using
the fact that the Brandt matrix of the Hecke correspondence $T_2$ has
column sums which are even. However, this evenness property fails in
general---for instance, when $p=23$ we have
$$
B(2) = \begin{bmatrix} 1 & 2 & 0 \\
1 & 1 & 1 \\
0 & 3 & 0
\end{bmatrix}
$$ which means we cannot expect $P_D(j_{23}(E))$ to be a perfect
square unless the number $j_{23}(E) \bmod 23$ differs from every
possible pair of supersingular $j_{23}$--invariants by quantities
having the same quadratic character mod $23$. This condition is
fulfilled by about one quarter of the curves satisfying the hypotheses
of Theorem~\ref{main}, and for these curves the proof of the theorem
goes through unchanged.

Even when $P_D(X)$ is not guaranteed to be a perfect square mod $p$,
empirical evidence indicates that the polynomial is sometimes a
perfect square anyway. For example, when $p=23$, a computer search up
to $\ell=400$ indicates that the primes $101$, $173$, and $317$ have
polynomials with square factorizations. It therefore seems reasonable
that classifying the square occurrences of $P_D(X)$ would lead to a
proof of the theorem in these cases.

\section{Acknowledgements}

This paper is based upon doctoral research performed under the
supervision of Noam Elkies, who first proposed to me the idea of
studying the curves $X_0^*(p)$ and provided valuable insight and
guidance during the time that I worked with him. The
work was supported under a National Defense Science and Engineering
Graduate Fellowship.


\begin{thebibliography}{99}
\bibitem{al} A. O. L. Atkin \& J. Lehner, Hecke operators on
$\Gamma_0(m)$, {\em Math. Ann.} {\bf 185} (1970), 134--160.
\bibitem{ant} Proceedings of the International Summer School on
Modular Functions of One Variable and Arithmetical Applications, RUCA,
University of Antwerp, Antwerp, July 17--August 3, 1972 (B. J. Birch \&
W. Kuyk eds.), Lecture Notes in Mathematics {\bf 476},
Springer--Verlag, Berlin, 1975. 
\bibitem{cox} D. Cox, {\em Primes of the form $x^2 + ny^2$}, John
Wiley \& Sons, Inc., New York, 1989.
\bibitem{deuring} M. Deuring, Die Typen der Multiplikatorenringe
elliptischer Funktionenk\"orper. {\em Abh. Math. Sem. Hansischen
Univ.} {\bf 14} (1941) 197--272.
\bibitem{e1} N. Elkies, The existence of infinitely many
supersingular primes for every elliptic curve over $\Q$. {\em
Invent. Math.} {\bf 89} (1987) 561--567.
\bibitem{e2} N. Elkies, Supersingular primes for elliptic
curves over real number fields. {\em Compositio Math.} {\bf 72} (1989)
165--172.
\bibitem{e3} N. Elkies, Elliptic and modular curves over finite
fields and related computational issues. {\em Computational
Perspectives on Number Theory: Proceedings of a Conference in Honor of
A.O.L. Atkin} (D.A. Buell and J.T. Teitelbaum, eds.), AMS/International
Press, 1998,
21--76.
\bibitem{gr1} B. Gross, Heegner Points on $X_0(N)$. {\em Modular Forms}
(R. A. Rankin, ed.), Ellis Horwood Company, Chichester, 1984, 87--105.
\bibitem{gz} B. Gross \& D. Zagier, On singular moduli.  {\em J. reine
angew.  Math.}  {\bf 355}  (1985), 191--220.
\bibitem{dj} D. Jao, {\em Supersingular Primes for Rational Points
on Modular Curves}, Harvard Doctoral Thesis, 2003.
\bibitem{kohel} D. Kohel, {\em Endomorphism rings of elliptic curves
  over finite fields}, Berkeley Doctoral Thesis, 1996.
\bibitem{lang} S. Lang, {\em Algebraic Number Theory, Second
Edition}, Springer--Verlag, 1994.
\bibitem{pizer} A. Pizer, {\em An Algorithm for computing modular
  forms on $\Gamma_0(N)$}. {\em J. Algebra} {\bf 64} (1980), 340--390.
\bibitem{aec} J. Silverman, {\em The Arithmetic of Elliptic
Curves}, Springer--Verlag, New York, 1986.
\end{thebibliography}
\end{document}